\documentclass[11pt]{article}
\usepackage{amsmath, amssymb, epsfig, graphicx}
\usepackage{bm}
\usepackage{mathrsfs}
\usepackage{color}
\usepackage{fullpage} 
\usepackage[small,bf]{caption}
\newtheorem{theorem}{Theorem}[section]

\newtheorem{lemma}[theorem]{Lemma}

\newtheorem{definition}[theorem]{Definition}

\numberwithin{equation}{section}

\newcommand{\bigO}{\mathrm{O}}

\DeclareMathOperator*{\argmin}{arg min}

\def \R {\mathbb{R}}

\def \P {\mathbb{P}}
\def \E {\mathbb{E}}
\newcommand{\e}{\mathrm{e}} 
\def \eps {\varepsilon}
\def \< {\langle}
\def \> {\rangle}

\def \endprf{\hfill {\vrule height6pt width6pt depth0pt}\medskip}

\usepackage{fullpage}

\begin{document}
\bibliographystyle{plain}

\pagestyle{plain}

\title{Compressed Sensing with Coherent and Redundant Dictionaries}
\author{Emmanuel J. Cand\`es$^1$\thanks{Corresponding author: Emmanuel J. Cand\`es. Email: candes@stanford.edu}, Yonina C. Eldar$^2$, Deanna Needell$^1$, and Paige Randall$^3$\\
  \vspace{-.1cm}\\
  $^1$Departments of Mathematics and Statistics, Stanford University, Stanford, CA 94305\\
    \vspace{-.3cm}\\
  $^2$Department of Electrical Engineering, Technion - Israel Institute of Technology, Haifa 32000\\
    \vspace{-.3cm}\\
  $^3$Center for Communications Research, Princeton, NJ 08540}
\date{May 2010; Revised October, 2010}

\maketitle

\vspace{-0.3in}

\begin{abstract}
  This article presents novel results concerning the recovery of
  signals from undersampled data in the common situation where such
  signals are not sparse in an orthonormal basis or incoherent
  dictionary, but in a truly redundant dictionary.  This work thus
  bridges a gap in the literature and shows not only that compressed
  sensing is viable in this context, but also that accurate recovery
  is possible via an $\ell_1$-analysis optimization problem.  We
  introduce a condition on the measurement/sensing matrix, which is a
  natural generalization of the now well-known restricted isometry
  property, and which guarantees accurate recovery of signals that are
  nearly sparse in (possibly) highly overcomplete and coherent
  dictionaries.  This condition imposes no incoherence restriction on
  the dictionary and our results may be the first of this kind.  We
  discuss practical examples and the implications of our results on
  those applications, and complement our study by demonstrating the
  potential of $\ell_1$-analysis for such problems.
\end{abstract}


\section{Introduction}

Compressed sensing is a new data acquisition theory based on the
discovery that one can exploit sparsity or compressibility when
acquiring signals of general interest, and that one can design
nonadaptive sampling techniques that condense the information in a
compressible signal into a small amount of
data~\cite{CRT06:Robust-Uncertainty,CT04:Near-Optimal,Don06:Compressed-Sensing}.
In a nutshell, reliable, nonadaptive data acquisition, with far fewer
measurements than traditionally assumed, is possible. By now,
applications of compressed sensing are abundant and range from imaging
and error correction to radar and remote sensing, see
\cite{IEEEMag:25-2,IEEEMag:24-4} and references therein.

In a nutshell, compressed sensing proposes acquiring a signal $x \in
\R^n$ by collecting $m$ linear measurements of the form $y_k = \langle
a_k, x\rangle + z_k$, $1 \le k \le m$, or in matrix notation,
\begin{equation}
\label{eq:model}
y = Ax + z;
\end{equation}
$A$ is an $m \times n$ sensing matrix with $m$ typically smaller than
$n$ by one or several orders of magnitude (indicating some significant
undersampling) and $z$ is an error term modeling measurement
errors. Sensing is nonadaptive in that $A$ does not depend on
$x$. Then the theory asserts that if the unknown signal $x$ is
reasonably sparse, or approximately sparse, it is possible to recover
$x$, under suitable conditions on the matrix $A$, by convex
programming: we simply find the solution to
\begin{equation}\tag{$L_1$}
  \min_{\tilde x \in \R^n} \|\tilde x\|_1 \quad\text{subject to}\quad \|A\tilde x - y\|_2 \leq \eps, 
\end{equation}
where $\|\cdot\|_2$ denotes the standard Euclidean norm, $\|x\|_1 =
\sum |x_i|$ is the $\ell_1$-norm and $\eps^2$ is a likely upper bound
on the noise power $\|z\|_2^2$. (There are other algorithmic
approaches to compressed sensing based on greedy algorithms such as
Orthogonal Matching
Pursuit~\cite{MZ93:Matching-Pursuits,TG07:Signal-Recovery}, Iterative
Thresholding~\cite{BD08:Iterative,FR07:Iterative}, Compressive
Sampling Matching Pursuit~\cite{NT08:Cosamp}, and many others.)

Quantitatively, a concrete example of a typical result in compressed
sensing compares the quality of the reconstruction from the data $y$
and the model \eqref{eq:model} with that available if one had an
oracle giving us perfect knowledge about the most significant entries
of the unknown signal $x$.  Define -- here and throughout -- by $x_s$
the vector consisting of the $s$ largest coefficients of $x \in \R^n$
in magnitude:
\begin{equation}
\label{eq:xs}
  x_s = \argmin_{\|\tilde x\|_0 \le s} \|x - \tilde x\|_2,
\end{equation}
where $\|x\|_0 = |\{i : x_i \neq 0\}|$. In words, $x_s$ is the best
$s$-sparse approximation to the vector $x$, where we shall say that a
vector is $s$-sparse if it has at most $s$ nonzero entries.  Put
differently, $x-x_s$ is the tail of the signal, consisting of the
smallest $n - s$ entries of $x$. In particular, if $x$ is $s$-sparse,
$x - x_s = 0$. Then with this in mind, one of the authors
\cite{Can08:Restricted-Isometry} improved on the work of Cand\`es,
Romberg and Tao~\cite{CRT06:Stable} and established that $(L_1)$
recovers a signal $\hat{x}$ obeying
\begin{equation}\label{eq:l1}
  \|\hat{x} - x\|_2 \leq C_0\frac{\|x-x_s\|_1}{\sqrt{s}} + C_1\eps, 
\end{equation}
provided that the {\em $2s$-restricted isometry constant} of $A$ obeys
$\delta_{2s} < \sqrt{2} - 1$.  The constants in this result have been
further improved, and it is now known to hold when $\delta_{2s} <
0.4652$~\cite{F10:Anote}, see also \cite{FL08:Sparsest}.  In short,
the recovery error from $(L_1)$ is proportional to the measurement
error and the tail of the signal.  This means that for compressible
signals, those whose coefficients obey a power law decay, the
approximation error is very small, and for exactly sparse signals it
completely vanishes.

The definition of restricted isometries first appeared in
\cite{CT05:Decoding} where it was shown to yield the error bound
\eqref{eq:l1} in the noiseless setting, i.~e.~when $\eps = 0$ and $z =
0$.
\begin{definition}
  For an $m\times n$ measurement matrix $A$, the
  $s$-\textit{restricted isometry constant} $\delta_s$ of $A$ is the
  smallest quantity such that
$$
(1-\delta_s)\|x\|_2^2 \leq \|Ax\|_2^2 \leq (1+\delta_s)\|x\|_2^2
$$
holds for all $s$-sparse signals $x$.
\end{definition}

With this, the condition underlying \eqref{eq:l1} is fairly natural
since it is interpreted as preventing sparse signals from being in the
nullspace of the sensing matrix $A$. Further, a matrix having a small
restricted isometry constant essentially means that every subset of
$s$ or fewer columns is approximately an orthonormal system.  It is
now well known that many types of random measurement matrices have
small restricted isometry
constants~\cite{CT04:Near-Optimal,MPJ06:Uniform,RV08:sparse,BDDW07:Johnson-Lindenstrauss}.
For example, matrices with Gaussian or Bernoulli entries have small
restricted isometry constants with very high probability whenever the
number of measurements $m$ is on the order of $s\log(n/s)$.  The fast
multiply matrix consisting of randomly chosen rows of the discrete
Fourier matrix also has small restricted isometry constants with very
high probability with $m$ on the order of $s (\log n)^4$.

\subsection{Motivation}
The techniques above hold for signals which are sparse in the standard
coordinate basis or sparse with respect to some other {\em orthonormal
  basis}.  However, there are numerous practical examples in which a
signal of interest is not sparse in an orthonormal basis.  More often
than not, sparsity is expressed not in terms of an orthonormal basis
but in terms of an \textit{overcomplete} dictionary.  This means that
our signal $f \in \R^n$ is now expressed as $f = Dx$ where $D \in
\R^{n \times d}$ is some overcomplete dictionary in which there are
possibly many more columns than rows.  The use of overcomplete
dictionaries is now widespread in signal processing and data analysis,
and we give two reasons why this is so. The first is that there may
not be any sparsifying orthonormal basis, as when the signal is
expressed using curvelets~\cite{CD02:New-Tight,CDDY05:Fast} or
time-frequency atoms as in the Gabor
representation~\cite{FS98:Gabor-Analysis}. In these cases and others,
no good orthobases are known to exist and researchers work with tight
frames. The second reason is that the research community has come to
appreciate and rely on the flexibility and convenience offered by
overcomplete representations. In linear inverse problems such as
deconvolution and tomography for example -- and even in straight
signal-denoising problems where $A$ is the identity matrix -- people
have found overcomplete representations to be extremely helpful in
reducing artifacts and mean squared error
(MSE)~\cite{SED04:Redundant,SFM07:uwd}.  It is only natural to expect
overcomplete representations to be equally helpful in compressed
sensing problems which, after all, are special inverse problems.
  
Although there are countless applications for which the signal of
interest is represented by some overcomplete dictionary, the
compressed sensing literature is lacking on the subject.  Consider the
simple case in which the sensing matrix $A$ has Gaussian (standard
normal) entries.  Then the matrix $AD$ relating the observed data with
the assumed (nearly) sparse coefficient sequence $x$ has independent
rows but each row is sampled from $\mathcal{N}(0, \Sigma)$, where
$\Sigma = D^*D$.  If $D$ is an orthonormal basis, then these entries
are just independent standard normal variables, but if $D$ is not
unitary then the entries are correlated, and $AD$ may no longer
satisfy the requirements imposed by traditional compressed sensing
assumptions.  In~\cite{RandThes} recovery results are obtained when the sensing
matrix $A$ is of the form $\Phi D^*$ where $\Phi$ satisfies the restricted
isometry property.  In this case the sampling matrix must depend on the
dictionary $D$ in which the signal is sparse.  We look for a universal result
which allows the sensing matrix to be independent from the signal and its
representation.  To be sure, we are not aware of any such results in the
literature guaranteeing good recovery properties when the columns may
be highly -- and even perfectly -- correlated.

Before continuing, it might be best to fix ideas to give some examples
of applications in which redundant dictionaries are of crucial
importance.
\begin{description}
\item [Oversampled DFT] The Discrete Fourier Transform (DFT) matrix is
  an $n \times n$ orthogonal matrix whose $k$th column is given by
\[
d_{k}(t) = \frac{1}{\sqrt{n}} e^{-2\pi i kt/n},
\]
with the convention that $0 \le t, k \le n-1$. Signals which are
sparse with respect to the DFT are only those which are superpositions
of sinusoids with frequencies appearing in the lattice of those in the
DFT.  In practice, we of course rarely encounter such signals.  To
account for this, one can consider the oversampled DFT in which the
sampled frequencies are taken over even smaller equally spaced
intervals, or at small intervals of varying lengths.  This leads to an
overcomplete frame whose columns may be highly correlated.

\item [Gabor frames] Recall that for a fixed function $g$ and positive
  time-frequency shift parameters $a$ and $b$, the $k$th column (where
  $k$ is the double index $k = (k_1, k_2)$) of the Gabor frame is
  given by
\begin{equation}\label{eq:gabor}
G_{k}(t) = g(t-k_2a) e^{2\pi ik_1bt}.
\end{equation}
Radar and sonar along with other imaging systems appear in many
engineering applications, and the goal is to recover pulse trains
given by
\[
f(t) = \sum_{j=1}^k \alpha_j
w\Big(\frac{t-t_j}{\sigma_j}\Big)\e^{i\omega_j t}. 
\]
Due to the time-frequency structure of these applications, Gabor
frames are widely used~\cite{Mal99:Wavelet-Tour}.  If one wishes to
recover pulse trains from compressive samples by using a Gabor
dictionary, standard results do not apply.

\item [Curvelet frames] Curvelets provide a multiscale decomposition of
  images, and have geometric features that set them apart from
  wavelets and the likes.  Conceptually, the curvelet transform is a
  multiscale pyramid with many directions and positions at each length
  scale, and needle-shaped elements at fine scales~\cite{CD02:New-Tight}.
  The transform gets its name from the fact that it approximates well
  the curved singularities in an image.  This transform has many
  properties of an orthonormal basis, but is overcomplete.  Written in
  matrix form $D$, it is a tight frame obeying the Parseval relations
$$
f = \sum_k \< f, d_k \> d_k \quad\text{and}\quad \|f\|_2^2 = \sum_k | \< f, d_k \> |^2,
$$
where we let $\{d_k\}$ denote the columns of $D$.  Although columns of
$D$ far apart from one another are very uncorrelated, columns close to
one another have high correlation.  Thus none of the results in
compressed sensing apply for signals represented in the curvelet
domain.

\item [Wavelet Frames] The undecimated wavelet transform (UWT) is a
  wavelet transform achieving translation invariance, a property that
  is missing in the discrete wavelet transform (DWT)~\cite{D89:uwt}.
  The UWT lacks the downsamplers and upsamplers in the DWT but
  upsamples the filter coefficients by a factor of $2^m$ at the
  $(m-1)$st level -- hence it is overcomplete.  Also, the
   Unitary Extension Principle of Ron and Shen~\cite{RS97:Affine} facilitates
   tight wavelet frame constructions for $L^2(\mathbb{R}^d)$ which
   may have many more wavelets than in the orthonormal case.
   This redundancy has
  been found to be helpful in image processing (see
  e.g.~\cite{SED04:Redundant}), and so one wishes for a recovery
  result allowing for significant redundancy and/or correlation.

\item [Concatenations] In many applications a signal may not be sparse
  in a single orthonormal basis, but instead is sparse over several
  orthonormal bases.  For example, a linear combination of spikes and
  sines will be sparse when using a concatenation of the coordinate
  and Fourier bases. 
  One also benefits by exploiting geometry and pointwise singularities
  in images by using combinations of tight frame coefficients such as
  curvelets, wavelets, and brushlets.  However, due to the correlation
  between the columns of these concatenated bases, current compressed
  sensing technology does not apply.
\end{description}

These and other applications strongly motivate the need for results
applicable when the dictionary is redundant and has correlations.
This state of affair, however, exposes a large gap in the literature
since current compressed sensing theory only applies when the
dictionary is an orthonormal basis, or when the dictionary is
extremely uncorrelated (see
e.g.~\cite{GMS03:Approximation-Functions,RSV08:Redundant,BCJ10:Why}).

\subsection{Do we really need incoherence?}

Current assumptions in the field of compressed sensing and sparse
signal recovery impose that the measurement matrix have uncorrelated
columns.  To be formal, one defines the \textit{coherence} of a matrix
$M$ as
$$
\mu(M) = \max_{j<k} \frac{|\< M_j, M_k \> |}{\|M_j\|_2 \|M_k\|_2},
$$
where $M_j$ and $M_k$ denote columns of $M$.  We say that a dictionary
is \textit{incoherent} if $\mu$ is small.  Standard results then
require that the measurement matrix satisfy a strict incoherence
property \cite{Tro04:Greed-Good,CP07:Near}, as even the RIP imposes
this.  If the dictionary $D$ is highly coherent, then the matrix $AD$
will also be coherent in general.

Coherence is in some sense a natural property in the compressed
sensing framework, for if two columns are closely correlated, it will
be impossible in general to distinguish whether the energy in the
signal comes from one or the other.\footnote{ Recall that when the
  dictionary $D$ is sufficiently incoherent, standard compressed
  sensing guarantees that we recover $x$ and thus $f=Dx$, provided $x$
  is $s$-sparse with $s$ sufficiently small. }  For example, imagine
that we are not undersampling and that $A$ is the identity so that we
observe $y = Dx$. Suppose the first two columns are identical, $d_1 =
d_2$. Then the measurement $d_1$ can be explained by the input vectors
$(1, 0, \ldots, 0)$ or $(0, 1, 0, \ldots, 0)$ or any convex
combination.  Thus there is no hope of reconstructing a unique sparse
signal $x$ from measurements $y = ADx$.  However, we are \textit{not}
interested in recovering the coefficient vector $x$, but rather the
actual signal $Dx$.  The large correlation between columns in $D$ now
does not impose a problem because although it makes it impossible to
tell apart coefficient vectors, this is not the goal. This simple
example suggests that perhaps coherence is not necessary.  If $D$ is
coherent, then we clearly cannot recover $x$ as in our example, but we
may certainly be able to recover the signal $f = Dx$ from measurements
$y=Af$ as we shall see next.

\begin{figure}[ht]
\begin{center}
  \includegraphics[scale=0.6]{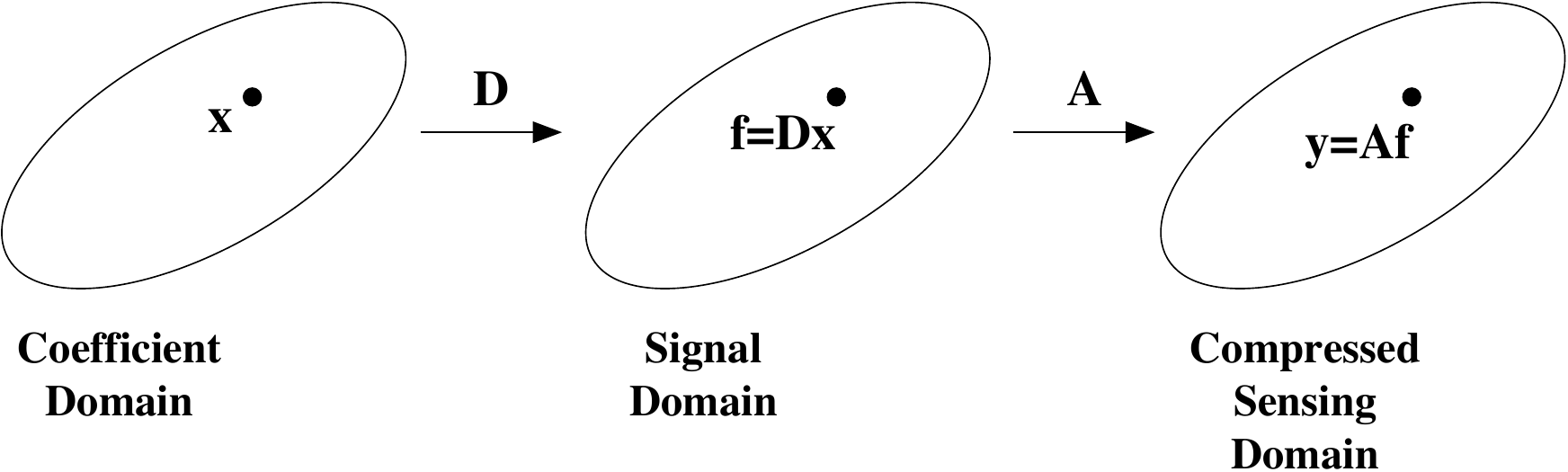}
  \caption{The compressed sensing process and its domains. This
    distinguishes the domains in which the measurements, signals, and
    coefficients reside}\label{fig:domains}
  \end{center}
\end{figure}

\subsection{Gaussian sensing matrices} 

To introduce our results, it might be best for pedagogical purposes to
discuss a concrete situation first, and we here assume that the
sensing matrix has iid Gaussian entries. In practice, signals are
never exactly sparse, and dictionaries are typically designed to make
$D^*f$ for some classes of $f$ as sparse as possible. Therefore, in
this paper, we propose a reconstruction from $y = A f + z$ by the
method of $\ell_1$-analysis:
\begin{equation}\tag{$P_1$}
  \hat{f} = \argmin_{\tilde f \in \R^n} \|D^* \tilde f\|_1 \quad\text{subject to}\quad \|A \tilde f - y\|_2 \leq \eps,
\end{equation}
where again $\eps$ is a likely upper bound on the noise level
$\|z\|_2$. Empirical studies have shown very promising results for the
$\ell_1$-analysis problem.  Its geometry has been
studied~\cite{EMR07:Analysis} as well as its applications to image
restoration~\cite{COS09:Split}.  However, there are no results in the
literature about its performance.

Our main result is that the solution to $(P_1)$ is very accurate
provided that $D^*f$ has rapidly decreasing coefficients. Our result
for the Gaussian case is below while the general theorem appears in
Section~\ref{sec:axiom}.
\begin{theorem}\label{thm:Gauss}
  Let $D$ be an arbitrary $n \times d$ tight frame and let $A$ be a
  $m\times n$ Gaussian matrix with $m$ on the order of $s\log(d/s)$.
  Then the solution $\hat{f}$ to $(P_1)$ obeys
$$
\|\hat{f} - f\|_2 \leq C_0\eps + C_1\frac{\|D^*f - (D^*f)_s\|_1}{\sqrt{s}},
$$
for some numerical constants $C_0$ and $C_1$, and where $(D^*f)_s$ is the vector consisting of the largest $s$ entries of $D^*f$ in magnitude as in~\eqref{eq:xs}.
\end{theorem}

We have assumed that $D$ is a tight frame although this is simply to
make the analysis easier and is of course not necessary.  Having said
this, our results proves not only that compressed sensing is viable
with highly coherent dictionaries, but also that the $\ell_1$-analysis
problem provably works in this setting.  We are not aware of any other
result of this kind. To be sure, other methods for redundant
dictionaries such as~\cite{RSV08:Redundant,BCJ10:Why,Tro04:Greed-Good}
force incoherence on the dictionary $D$ so that the matrix $AD$
conforms to standard compressed sensing results.  The method in~\cite{RandThes}
requires that the sensing matrix depend on the dictionary.  These are drastically
different from the setting here, where we impose no such 
properties on the dictionary.  We point out that our result holds even
when the coherence of the dictionary $D$ is maximal, meaning two
columns are completely correlated.  Finally, we also note that the
dependence on the noise level is optimal and that the tail bound in
the error is analogous to previous bounds in the non-redundant case
such as~\eqref{eq:l1}.

\subsection{Implications}

As we mentioned, the dependence on the noise in the error given by
Theorem~\ref{thm:Gauss} is optimal, and so we need only discuss how
the second term affects the estimation error.  This term will of
course be negligible when the norm of the tail, $D^*f - (D^*f)_s$, is
small.  Hence, the result says that for any dictionary, signals $f$
such that $D^*f$ decays rapidly can be approximately reconstructed
using $\ell_1$-analysis from just a few random measurements.  This is
exactly the case for many dictionaries used in practice and many
classes of signals as discussed earlier. As a side remark, one can
also guarantee rapid decay of $D^*f$ (we assume the signal expansion
$f = Dx$) when $D^*D$ is well behaved and the coefficient vector $x$
is nearly sparse. To see why this is true, suppose $D$ is a tight
frame so that $D^* f = D^* D x$. A norm commonly used to quantify
sparsity is the quasi $p$-norm with $p \le 1$ defined via $\|x\|_p^p
= \sum_i |x_i|^p$ (sparser signals with unit $2$-norm have
smaller $p$-norms). Now a simple calculation shows that
\[
\|D^* f\|_p \le \Bigl[ \max_j \sum_i |(D^*D)_{ij}|^p \Bigr]^{1/p} \,
\|x\|_p. 
\]
In words, if the columns of the Gram matrix are reasonably sparse and
if $f$ happens to have a sparse expansion, then the frame coefficient
sequence $D^* f$ is also sparse. All the transforms discussed above,
namely, the Gabor, curvelet, wavelet frame, oversampled Fourier
transform all have nearly diagonal Gram matrices -- and thus, sparse
columns.

We now turn to the implications of our result to the applications we
have already mentioned, and instantiate the theorem in the noiseless
case due to the optimality of the noise level in the error.

\begin{description}

\item[Multitone signals] To recover multitone signals, we use an
  oversampled DFT, which is not orthonormal and may have very large
  coherence. However, since each ``off-grid'' tone has a rapidly
  decaying expansion, $D^*f$ will have rapidly decaying
  coefficients.\footnote{In practice, one smoothly localizes the data
    to a time interval by means of a nice window $w$ to eliminate
    effects having to do with a lack of periodicity. One can then
    think of the trigonometric exponentials as smoothly vanishing at
    both ends of the time interval under study.} Thus our result
  implies that the recovery error is negligible when the number of
  measurements is about the number of tones times a log factor.

\item[Radar] For radar and sonar applications using Gabor
  dictionaries, our result similarly implies a negligible error.
  Indeed, with notation as in~\eqref{eq:gabor}, one sees that the
  sequence $\{\< w(t)\e^{i\omega t}, G_k(t)\rangle\}_k$ decays quickly
  (each pulse has a rapidly decaying expansion).  Therefore, our
  result implies a negligible error when the number of measurements is
  roughly the number of pulses in the pulse train, up
  to a log factor.

\item[Images] Roughly speaking, the curvelet coefficient sequence of
  an arbitrary image, which is discontinuous along piecewise-$C^2$
  edges but is otherwise smooth, decays like $k^{-3/2}$ -- up to a log
  factor -- when arranged in a decreasing order of magnitude. Hence,
  our theorem asserts that one can get an $\ell_2$ error of about
  $s^{-1}$ from about $s \log n$ random samples of $f$. This is
  interesting since this is the approximation error one would get by
  bandlimiting the image to a spatial frequency about equal to $s$ or,
  equivalently, by observing the first $s^2$ Fourier coefficients of
  $f$. So even though we do not know where the edges are, this says
  that one can sense an image nonadaptively $m$ times, and get a
  recovery error which is as good as that one would traditionally get
  by taking a much larger number -- about $m^2$ -- of samples. This is
  a considerable gain. Of course, similar implications hold when the
  undecimated wavelet transform of those images under study decay
  rapidly.

\item[Concatenations] When working with signals which are sparse over
  several orthonormal bases, it is natural to use a dictionary $D$
  consisting of a concatenation of these bases.  For example, consider
  the dictionary $D$ which is the concatenation of the identity and
  the Fourier basis (ignoring normalizations for now).  Then $D^*D$ is
  made up of four blocks, two of which are the identity and two of
  which are the DFT, and does not have sparse columns.  Then even when
  $f$ is sparse in $D$, the coefficients of $D^*f$ may be spread.  If
  this is the case, then the theorem does not provide a good error
  bound.  This should not be a surprise however, for if $D^*f$ is not
  close to a sparse signal, then we do not expect $f$ to be the
  minimizer of the $\ell_1$-norm in $(P_1)$.  In this case,
  $\ell_1$-analysis is simply not the right method to use.
\end{description}

To summarize, we see that in the majority of the applications, our
theorem yields good recovery results.  As seen in the last example,
the $\ell_1$-analysis method only makes sense when $D^*f$ has quickly
decaying coefficients, which may not be the case for concatenations of
orthonormal bases.  However, this is not always the case, as we see in
the following.

\textbf{An easy example.}  As above, let $D$ be the $n \times 2n$
dictionary consisting of a concatenation of the identity and the DFT,
normalized to ensure $D$ is a tight frame (below, $F$ is the DFT
normalized to be an isometry):
\[
D = \frac{1}{\sqrt{2}}[ I\quad F ].
\]
We wish to create a sparse signal that uses linearly dependent columns
for which there is no local isometry. Assume that $n$ is a perfect
square and consider the Dirac comb
\[
f(t) = \sum_{j=1}^{\sqrt{n}} \delta(t - j\sqrt{n}),
\]
which is a superposition of spikes spread $\sqrt{n}$ apart.  
Thus our signal is a
sparse linear combination of spikes and sines, something that by the
last example alone we would not expect to be able to recover.  However,
$D^*f = [f \quad f]/\sqrt{2}$ is exactly sparse implying that $\|D^*f
- (D^*f)_s\|_1 = 0$ when $s > 2\sqrt{n}$. Thus our result shows that
$\ell_1$-analysis can exactly recover the Dirac comb consisting of
spikes and sines from just a few general linear functionals. 

\subsection{Axiomization}\label{sec:axiom}
We now turn to the generalization of the above result, and give broader conditions about the sensing matrix under which the recovery algorithm performs well.   We will impose a natural property on the measurement matrix, analogous to the restricted isometry property.

\begin{definition}[D-RIP]
  Let $\Sigma_s$ be the union of all subspaces spanned by all subsets
  of $s$ columns of $D$.  We say that the measurement matrix $A$ obeys
  the \textit{restricted isometry property adapted to $D$}
  (abbreviated D-RIP) with constant $\delta_s$ if
$$
(1-\delta_s)\|v\|_2^2 \leq \|Av\|_2^2 \leq (1+\delta_s)\|v\|_2^2
$$
holds for all $v \in \Sigma_s$. 
\end{definition}  

We point out that $\Sigma_s$ is just the image under $D$ of all
$s$-sparse vectors.  Thus the D-RIP is a natural extension to the
standard RIP.  We will see easily that Gaussian matrices and other
random compressed sensing matrices satisfy the D-RIP.  In fact any $m
\times n$ matrix $A$ obeying for fixed $v\in\R^n$,
\begin{equation}\label{eq:cond}
\P\big( (1-\delta)\|v\|_2^2 \leq \|Av\|_2^2 \leq (1+\delta)\|v\|_2^2 \big) \leq C\e^{-\gamma m}
\end{equation}
($\gamma$ is an arbitrary positive numerical constant) will satisfy
the D-RIP with overwhelming probability, provided that $m \gtrsim
s\log (d/s)$.  This can be seen by a standard covering argument (see
e.g.~the proof of Lemma 2.1 in~\cite{RSV08:Redundant}).  Many types of
random matrices satisfy~\eqref{eq:cond}.  It is now well known that
matrices with Gaussian, subgaussian, or Bernoulli entries
satisfy~\eqref{eq:cond} with number of measurements $m$ on the order
of $s\log (d/s)$ (see e.g.~\cite{BDDW07:Johnson-Lindenstrauss}).  It
has also been shown~\cite{MPJ06:Uniform} that if the rows of $A$ are
independent (scaled) copies of an isotropic $\psi_2$ vector, then $A$
also satisfies~\eqref{eq:cond}.  Recall that an isotropic $\psi_2$
vector $a$ is one that satisfies for all $v$,
$$
\E |\< a, v \> |^2 = \|v\|^2 \quad\text{and}\quad \inf \{t : \E\exp(\< a, v\> ^2 / t^2) \leq 2\} \leq \alpha\|v\|_2,
$$
for some constant $\alpha$.  See~\cite{MPJ06:Uniform} for further
details.  Finally, it is clear that if $A$ is any of the above random
matrices then for any fixed unitary matrix $U$, the matrix $AU$ will
also satisfy the condition.

The D-RIP can also be analyzed via the Johnson-Lindenstrauss lemma (see e.g.~\cite{HV09:johnson,AC09:fast}).  There are many results that show certain types of matrices satisfy this lemma, and these would then satisfy the D-RIP via~\eqref{eq:cond}.  Subsequent to our submission of this manuscript, Ward and Krahmer showed that randomizing the column signs of any matrix that satisfies the standard RIP yields a matrix which satisfies the Johnson-Lindenstrauss lemma~\cite{KW10:NewAnd}.  Therefore, nearly all random matrix constructions which satisfy standard RIP compressed sensing requirements will also satisfy the D-RIP.  A particularly important consequence is that because the randomly subsampled Fourier matrix is known to satisfy the RIP, this matrix along with a random sign matrix will thus satisfy D-RIP.  This gives a fast transform which satisfies the D-RIP.  See Section~\ref{sec:Fast} for more discussion.

We are now prepared to state our main result.

\begin{theorem}\label{thm:main}
Let $D$ be an arbitrary tight frame and let $A$ be a measurement matrix satisfying D-RIP with $\delta_{2s} < 0.08$.  Then the solution $\hat{f}$ to $(P_1)$ satisfies
$$
\|\hat{f} - f\|_2 \leq C_0\eps + C_1\frac{\|D^*f - (D^*f)_s\|_1}{\sqrt{s}},
$$
where the constants $C_0$ and $C_1$ may only depend on $\delta_{2s}$.
\end{theorem}

{\em Remarks.} We actually prove that the theorem holds under the
weaker condition $\delta_{7s} \leq 0.6$, however we have not tried to
optimize the dependence on the values of the restricted isometry
constants; refinements analagous to those in the compressed sensing
literature are likely to improve the condition.  Further, we note that
since Gaussian matrices with $m$ on the order of $s\log(d/s)$ obey the
D-RIP, Theorem~\ref{thm:Gauss} is a special case of
Theorem~\ref{thm:main}.

\subsection{Organization} 
The rest of the paper is organized as follows.  In
Section~\ref{sec:proof} we prove our main result,
Theorem~\ref{thm:main}.  Section~\ref{sec:nums} contains numerical
studies highlighting the impact of our main result on some of the
applications previously mentioned.  In Section~\ref{sec:disc} we
discuss further the implications of our result along with its
advantages and challenges.  We compare it to other methods proposed in
the literature and suggest an additional method to overcome some
impediments.  

\section{Proof of Main Result}\label{sec:proof}
We now begin the proof of Theorem~\ref{thm:main}, which is inspired by
that in~\cite{CRT06:Stable}.  The new challenge here is that although
we can still take advantage of sparsity, the vector possessing the sparse
property is not being multiplied by something that satisfies the RIP, as in 
the standard compressed sensing case.  Rather than bounding
the tail of $f-\hat{f}$ by its largest coefficients as in~\cite{CRT06:Stable}, we bound a portion of $D^*h$ in 
an analagous way.  We then utilize the D-RIP and the fact that $D$ is a tight
frame to bound the error, $\|f-\hat{f}\|_2$.

Let $f$ and $\hat{f}$ be as in the
theorem, and let $T_0$ denote the set of the largest $s$ coefficients
of $D^*f$ in magnitude.  We will denote by $D_T$ the matrix $D$ restricted to the
columns indexed by $T$, and write $D_T^*$ to mean $(D_T)^*$.   With $h = f - \hat{f}$, our goal is to bound
the norm of $h$.  We will do this in a sequence of short lemmas.  The first
is a simple consequence of the fact that $\hat{f}$ is the minimizer.

\begin{lemma}[Cone Constraint]\label{cone}
The vector $D^*h$ obeys the following cone constraint,
\begin{equation*}
\|D_{T_0^c}^*h\|_1 \leq 2\|D_{T_0^c}^*f\|_1 + \|D_{T_0}^*h\|_1.
\end{equation*}
\end{lemma}

\textit{Proof. }
Since both $f$ and $\hat{f}$ are feasible but
$\hat{f}$ is the minimizer, we must have $\|D^*\hat{f}\|_1 \leq
\|D^*f\|_1$.   We then
have that
\begin{align*}
\|D_{T_0}^*f\|_1 + \|D_{T_0^c}^*f\|_1 = \|D^*f\|_1 
&\geq \|D^*\hat{f}\|_1 \\
&= \|D^*f - D^*h\|_1 \\
&\geq \|D_{T_0}^*f\|_1 - \|D_{T_0}^*h\|_1 - \|D_{T_0^c}^*f\|_1 + \|D_{T_0^c}^*h\|_1.
\end{align*}
This implies the desired cone constraint.
\endprf

We next divide the coordinates $T_0^c$ into sets of size $M$ (to be
chosen later) in order of decreasing magnitude of $D_{T_0^c}^*h$.
Call these sets $T_1, T_2, \ldots$, and for simplicity of notation set
$T_{01} = T_0 \cup T_1$.  We then bound the tail of $D^*h$.

\begin{lemma}[Bounding the tail]\label{thirteen}
Setting $\rho = s/M$ and $\eta = 2\|D_{T_0^c}^*f\|_1 / \sqrt{s}$, we have the following bound,
\begin{equation*}
\sum_{j\geq 2} \|D^*_{T_j}h\|_2 \leq \sqrt{\rho}(\|D_{T_0}^*h\|_2 + \eta).
\end{equation*}
\end{lemma}
\textit{Proof. }
 By construction of the sets $T_j$, we have that each
coefficient of $D_{T_{j+1}}^*h$, written $|D_{T_{j+1}}^*h|_{(k)}$, is
at most the average of those on $T_j$:
$$
|D_{T_{j+1}}^*h|_{(k)} \leq \|D_{T_j}^*h\|_1 / M.
$$
Squaring these terms and summing yields
$$
\|D_{T_{j+1}}^*h\|_2^2 \leq \|D_{T_j}^*h\|_1^2 / M.
$$
This along with the cone constraint in Lemma~\ref{cone} gives
\[
\sum_{j\geq 2} \|D^*_{T_j}h\|_2 \leq \sum_{j\geq 1} \|D_{T_j}^*h\|_1 /
\sqrt{M} = \|D_{T_0^c}^*h\|_1 / \sqrt{M}.
\]
With $\rho = s/M$ and $\eta = 2\|D_{T_0^c}^*f\|_1 / \sqrt{s}$,
it follows from Lemma~\ref{cone} and the
Cauchy-Schwarz inequality that
\begin{equation*}
\sum_{j\geq 2} \|D^*_{T_j}h\|_2 \leq \sqrt{\rho}(\|D_{T_0}^*h\|_2 + \eta),
\end{equation*}
as desired.
\endprf

Next we observe that by the feasibility of $\hat{f}$, $Ah$ must be
small.  
\begin{lemma}[Tube Constraint]\label{Ah}
The vector $Ah$ satisfies the following,
$$
\|Ah\|_2 \leq 2\varepsilon.
$$
\end{lemma}
\textit{Proof. }
Since $\hat{f}$ is feasible, we have
\begin{equation*}
\|Ah\|_2 = \|Af - A\hat{f}\|_2 \leq \|Af - y\|_2 + \|A\hat{f} - y\|_2 \leq \varepsilon + \varepsilon = 2\varepsilon.
\end{equation*}
\endprf

We will now need the following result which utilizes the fact that $D$ satisfies the D-RIP.

\begin{lemma}[Consequence of D-RIP]\label{sixteen}
The following inequality holds,
\begin{equation*}
\sqrt{1 - \delta_{s+M}}\|D_{T_{01}}D_{T_{01}}^*h\|_2 - \sqrt{\rho(1 + \delta_M)}(\|h\|_2 + \eta) \leq 2\varepsilon.
\end{equation*}
\end{lemma}
\textit{Proof. }
Since $D$ is a tight frame, $DD^*$ is the identity, and this along with
the D-RIP and~Lemma~\ref{thirteen} then imply the following:
\begin{align*}
2\varepsilon &\geq \|Ah\|_2 = \|ADD^*h\|_2  \\
&\geq \|AD_{T_{01}}D_{T_{01}}^*h\|_2 - \sum_{j\geq 2} \|AD_{T_{j}}D_{T_{j}}^*h\|_2 \\
&\geq \sqrt{1 - \delta_{s+M}}\|D_{T_{01}}D_{T_{01}}^*h\|_2 - \sqrt{1 + \delta_M}\sum_{j\geq 2} \|D_{T_{j}}D_{T_{j}}^*h\|_2 \\
&\geq \sqrt{1 - \delta_{s+M}}\|D_{T_{01}}D_{T_{01}}^*h\|_2 - \sqrt{\rho(1 + \delta_M)}(\|D_{T_0}^*h\|_2 + \eta).
\end{align*}
 
Since we also have $\|D_{T_0}^*h\|_2 \leq \|h\|_2$, this yields the desired result.
\endprf

We now translate these bounds to the bound of the actual error, $\|h\|_2$.

\begin{lemma}[Bounding the error]\label{side}
The error vector $h$ has norm that satisfies,
$$
\|h\|_2^2 \leq \|h\|_2\|D_{T_{01}}D_{T_{01}}^*h \|_2 + \rho(\|D_{T_0}^*h\|_2 + \eta)^2,
$$
\end{lemma}
 
\textit{Proof. }
Since $D^*$ is an isometry, we have
\begin{align*}
\|h\|_2^2 &= \|D^*h\|_2^2 = \|D_{T_{01}}^*h\|_2^2 + \|D_{T_{01}^c}^*h\|_2^2\notag \\
&= \langle h, D_{T_{01}}D_{T_{01}}^*h \rangle + \|D_{T_{01}^c}^*h\|_2^2\notag \\
&\leq \|h\|_2\|D_{T_{01}}D_{T_{01}}^*h \|_2 + \|D_{T_{01}^c}^*h\|_2^2\notag \\
&\leq \|h\|_2\|D_{T_{01}}D_{T_{01}}^*h \|_2 + \rho(\|D_{T_0}^*h\|_2 + \eta)^2,
\end{align*}
where the last inequality follows from Lemma~\ref{thirteen}.  
\endprf


We next
observe an elementary fact that will be useful. The proof is omitted.
\begin{lemma}\label{lem}
  For any values $u$, $v$ and $c>0$, we have
$$
uv \leq \frac{cu^2}{2} + \frac{v^2}{2c}.
$$
\end{lemma}

We may now conclude the proof of Theorem~\ref{thm:main}. First we employ Lemma~\ref{lem} twice to the inequality given
by Lemma~\ref{side}
(with constants $c_1$, $c_2$ to be chosen later) and the bound
$\|D_{T_0}^*h\|_2 \leq \|h\|_2$ to get

\begin{align*}
  \|h\|_2^2 &\leq \frac{c_1\|h\|_2^2}{2} + \frac{\|D_{T_{01}}D_{T_{01}}^*h \|_2^2}{2c_1} + \rho(\|h\|_2 + \eta)^2\\
  & = \frac{c_1\|h\|_2^2}{2} + \frac{\|D_{T_{01}}D_{T_{01}}^*h \|_2^2}{2c_1} + \rho\|h\|_2^2 + 2\rho\eta\|h\|_2 + \rho\eta^2\\
  & \leq \frac{c_1\|h\|_2^2}{2} + \frac{\|D_{T_{01}}D_{T_{01}}^*h
    \|_2^2}{2c_1} + \rho\|h\|_2^2 + 2\rho\Big(\frac{c_2\|h\|_2^2}{2} +
  \frac{\eta^2}{2c_2} \Big) + \rho\eta^2.
\end{align*}
Simplifying, this yields
$$
\Big(1 - \frac{c_1}{2} - \rho - \rho c_2 \Big)\|h\|_2^2 \leq \frac{1}{2c_1}\|D_{T_{01}}D_{T_{01}}^*h \|_2^2 + \Big(\frac{\rho}{c_2} + \rho\Big)\eta^2.
$$

Using the fact that $\sqrt{u^2 + v^2} \leq u + v$ for $u,v\geq0$, we can further simply to get our desired lower bound,
\begin{equation}\label{lb}
\|D_{T_{01}}D_{T_{01}}^*h \|_2 \geq \|h\|_2\sqrt{2c_1\Big(1 - (\frac{c_1}{2}+\rho+\rho c_2) \Big)} - \eta\sqrt{2c_1\big(\frac{\rho}{c_2}+\rho\big)}.
\end{equation}

Combining~\eqref{lb} with Lemma~\ref{sixteen} implies
$$
2\varepsilon \geq K_1\|h\|_2 - K_2\eta,
$$
where
\begin{align*}
K_1 &= \sqrt{2c_1(1-\delta_{s+M})\Big(1 - (\frac{c_1}{2}+\rho+\rho c_2) \Big)} - \sqrt{\rho(1+\delta_{M})}, \quad \text{and}\\
K_2 &= \sqrt{2c_1(1-\delta_{s+M})(\rho/c_2 + \rho )} - \sqrt{\rho(1+\delta_{M})}.
\end{align*}

It only remains to choose the parameters $c_1$, $c_2$, and $M$ so that $K_1$ is positive.  We choose $c_1$=1, $M=6s$, and take $c_2$ arbitrarily small so that $K_1$ is positive when $\delta_{7s} \leq 0.6$.  Tighter restrictions on $\delta_{7s}$ will of course force the constants in the error bound to be smaller.  For example, if we set $c_1 = 1/2$, $c_2 = 1/10$, and choose $M=6s$, we have that whenever $\delta_{7s} \leq 1/2$ that $(P_1)$ reconstructs $\hat{f}$ satisfying
$$
\|f-\hat{f}\|_2 \leq 62\varepsilon + 30\frac{\|D_{T_0^c}^*f\|_1}{\sqrt{s}}.
$$
Note that if $\delta_{7s}$ is even a little smaller, say $\delta_{7s} \leq 1/4$, the constants in the theorem are just $C_1 = 10.3$ and $C_2 = 7.33$.  Note further that by Corollary 3.4 of~\cite{NT08:Cosamp}, $\delta_{7s} \leq 0.6$ is satisfied whenever $\delta_{2s} \leq 0.08$.  This completes the proof.
\endprf

\section{Numerical Results}\label{sec:nums}

We now present some numerical experiments illustrating the
effectiveness of recovery via $\ell_1$-analysis and also compare the
method to other alternatives.  Our results confirm that in practice,
$\ell_1$-analysis reconstructs signals represented in truly redundant
dictionaries, and that this recovery is robust with respect to noise.

\begin{figure}[ht]
\begin{center}
  \includegraphics[scale=0.4]{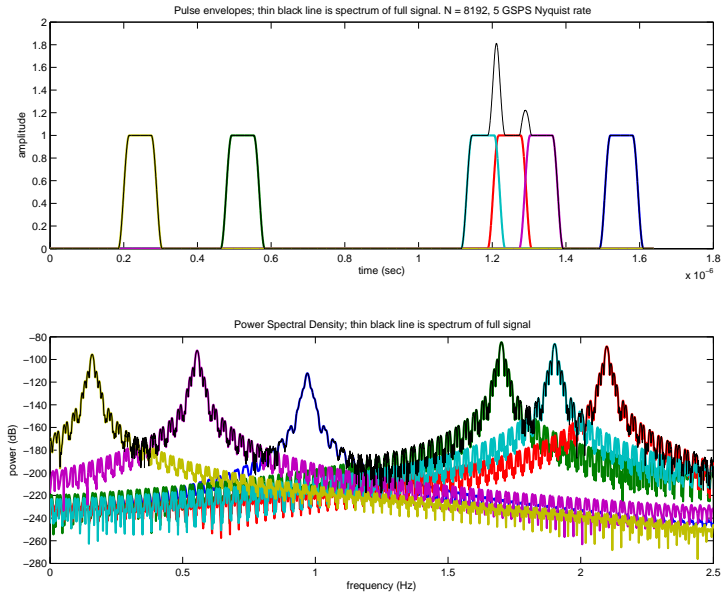}
\end{center}
\caption{Input signal in the time and frequency domains.  The signal
  of interest is a superposition of 6 radar pulses, each of which
  being about 200 ns long, and with frequency carriers distributed
  between 50 MHz and 2.5 GHz (top plot).  As can be seen, three of
  these pulses overlap in the time domain.}
\label{fig:input}
\end{figure}

In these experiments, we test the performance on a simulated
real-world signal from the field of radar detection. The test input is
a superposition of six radar pulses. Each pulse has a duration of
about 200 ns, and each pulse envelope is trapezoidal, with a 20 ns
rise and fall time, see Figure~\ref{fig:input}. For each pulse, the
carrier frequency is chosen uniformly at random from the range 50 MHz
to 2.5 GHz. The Nyquist interval for such signals is thus 0.2
ns. Lastly, the arrival times are distributed at random in a time
interval ranging from $t = 0 \text{ s}$ to $t \approx 1.64 \,
\mu\text{s}$; that is, the time interval under study contains $n =
8192$ Nyquist intervals. We acquire this signal by taking $400$
measurements only, so that the sensing matrix $A$ is a Gaussian matrix
with $400$ rows. The dictionary $D$ is a Gabor dictionary with
Gaussian windows, oversampled by a factor of about $60$ so that $d \approx 
60 \times 8,192 = 491,520$. The main comment about this setup is that
the signal of interest is not exactly sparse in $D$ since each pulse
envelope is not Gaussian (the columns of $D$ are pulses with Gaussian
shapes) and since both the frequencies and arrival times are sampled
from a continuous grid (and thus do not match those in the
dictionary).

\begin{figure}[ht]
\begin{center}
  \includegraphics[scale=0.4]{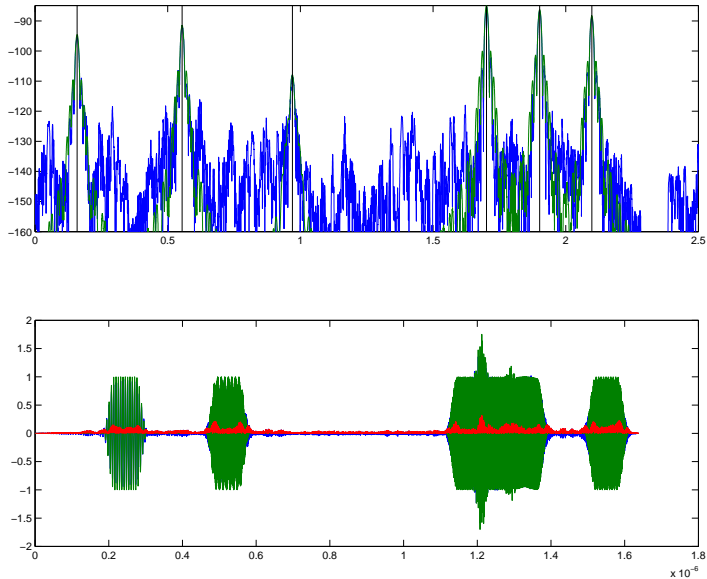}
\end{center}
\caption{Recovery in both the time (below) and frequency (above)
  domains by $\ell_1$-analysis. Blue denotes the recovered signal,
  green the actual signal, and red the difference between the two.}
\label{fig:timefreq}
\end{figure}

Figure~\ref{fig:timefreq} shows the recovery (without noise) by
$\ell_1$-analysis in both the time and frequency domains.  In the time
domain we see (in red) that the difference between the actual signal
and the recovered signal is small, as we do in the frequencey domain
as well.  These pulses together with the carrier frequencies are well
recovered from a very small set of measurements.

\begin{figure}[ht] 
  \begin{center}
    \includegraphics[scale=0.4]{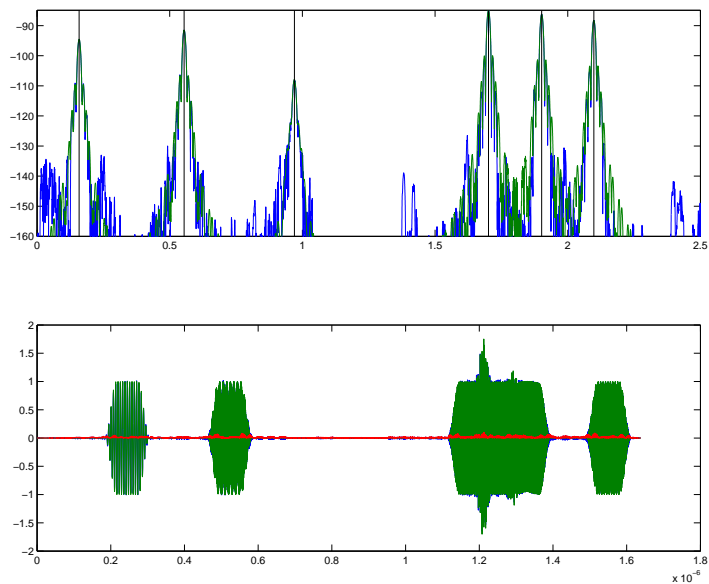}
\end{center}
\caption{Recovery in both the time (below) and frequency (above)
  domains by $\ell_1$-analysis after one reweighted iteration. Blue
  denotes the recovered signal, green the actual signal, and red the
  difference between the two. The RMSE is less than a third of that in
  Figure~\ref{fig:timefreq2}}\label{fig:timefreq2}
\end{figure}

In practice, reweighting the $\ell_1$ norm often offers superior
results. We use the \textit{reweighted} $\ell_1$-analysis method,
which solves several sequential weighted $\ell_1$-minimization
problems, each using weights computed from the solution of the
previous problem~\cite{CWB08:Reweighted}.  This procedure has been
observed to be very effective in reducing the number of measurements
needed for recovery, and outperforms standard $\ell_1$-minimization in
many situations (see
e.g.~\cite{CWB08:Reweighted},~\cite{KXAH:2step},~\cite{N09:rw}).
Figure~\ref{fig:timefreq2} shows reconstruction results after just one
reweighting iteration; the root-mean squared error (RMSE) is
significantly reduced, by a factor between 3 and 4.

\begin{figure}[ht]
\begin{center}
\begin{tabular}{cc}
\includegraphics[scale=0.25]{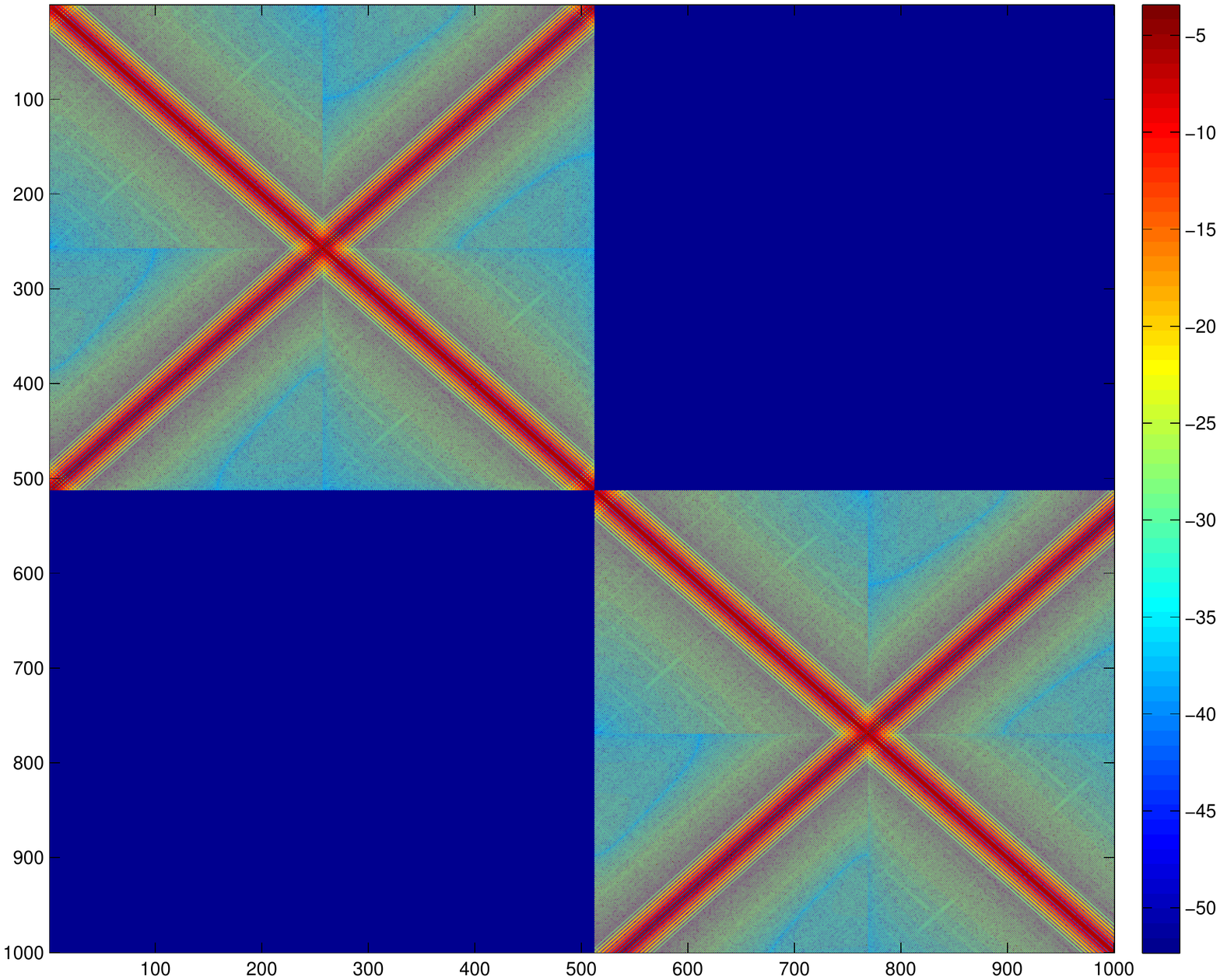}  & 
\includegraphics[scale=0.25]{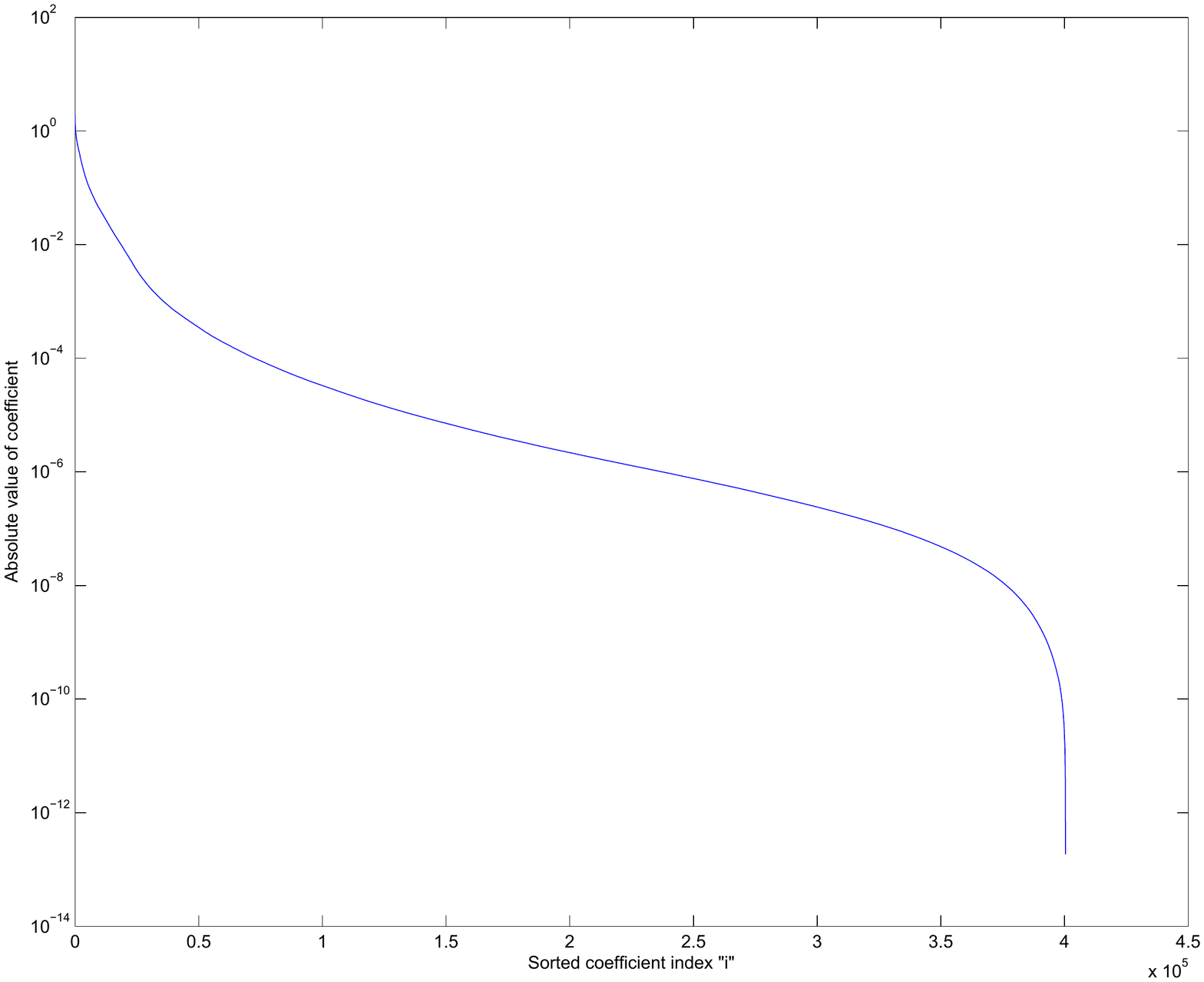}
\end{tabular}
\end{center}
\caption{Portion of the matrix $D^*D$, in $\log$-scale (left). Sorted
  analysis coefficients (in absolute value) of the signal from
  Figure~\ref{fig:input} (right).}\label{fig:notONB}
\end{figure}

Because $D$ is massively overcomplete, the Gram matrix $D^*D$ is not
diagonal. Figure~\ref{fig:notONB} depicts part of the Gram matrix
$D^*D$ for this dictionary, and shows that this matrix is ``thick''
off of the diagonal.  We can observe visually that the dictionary $D$
is not an orthogonal system or even a matrix with low coherence, and
that columns of this dictionary are indeed highly correlated. Having
said this, the second plot in Figure~\ref{fig:notONB} shows the rapid
decay of the sequence $D^*f$ where $f$ is the signal in
Figure~\ref{fig:input}.  

\begin{figure}[ht] 
  \begin{center}
 \includegraphics[scale=0.3]{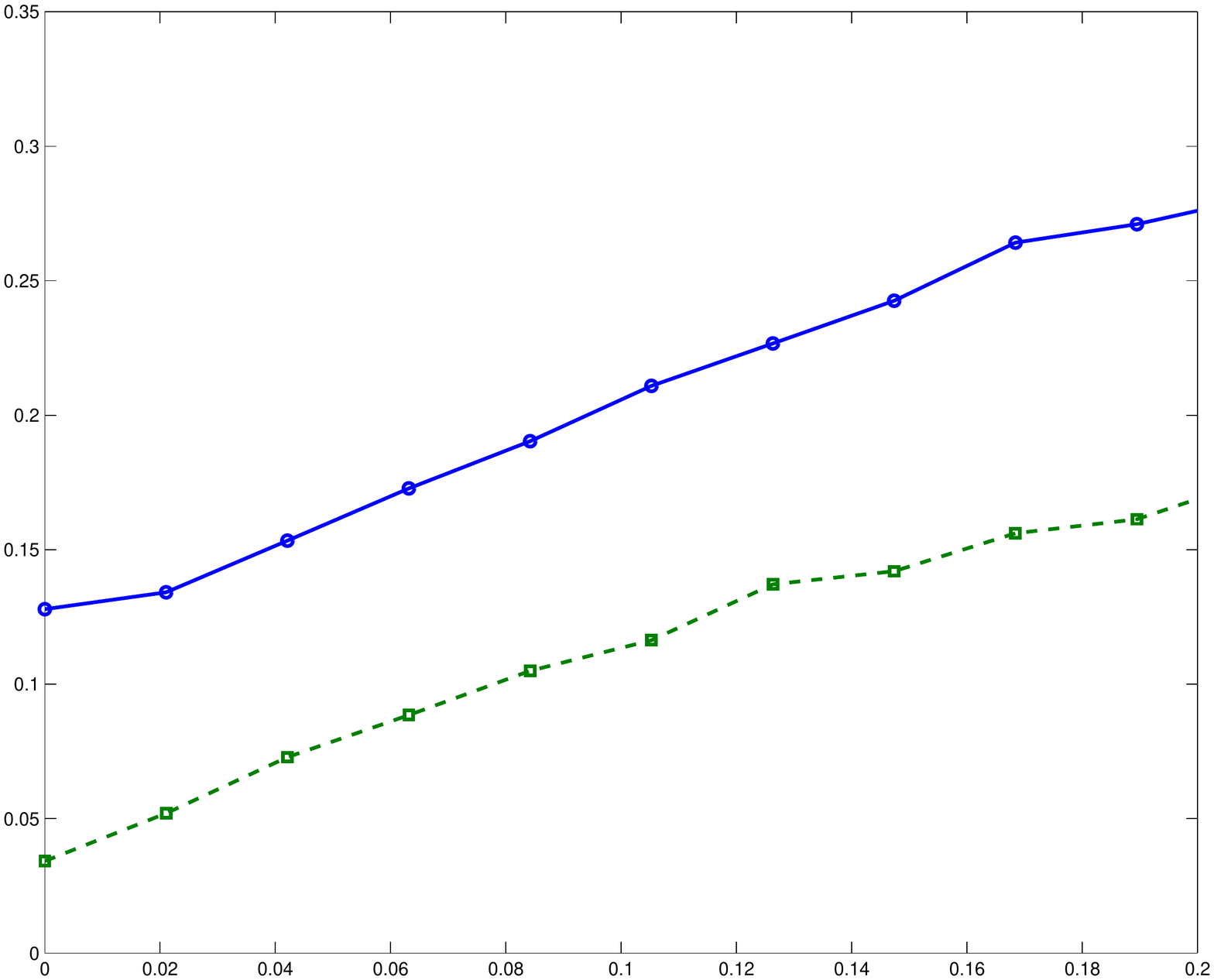}
\end{center}
\caption{Relative recovery error of $\ell_1$-analysis as a function of
  the (normalized) noise level, averaged over $5$ trials. The solid
  line denotes standard $\ell_1$-analysis, and the dashed line denotes
  $\ell_1$-analysis with $3$ reweighted iterations. The $x$-axis is
  the relative noise level $\sqrt{m} \sigma/\|Af\|_2$ while the
  $y$-axis is the relative error $\|\hat f -
  f\|_2/\|f\|_2$.}\label{fig:noise}
\end{figure}

Our next simulation studies the robustness of $\ell_1$-analysis with
respect to noise in the measurements $y = Af + z$, where $z$ is a
white noise sequence with standard deviation
$\sigma$. Figure~\ref{fig:noise} shows the recovery error as a
function of the noise level.  As expected, the relationship is linear,
and this simulation shows that the constants in Theorem~\ref{thm:main}
seem to be quite small.  This plot also shows the recovery error with
respect to noise using a \textit{reweighted} $\ell_1$-analysis;
reweighting also improves performance of $\ell_1$-analysis, as is seen
in Figure~\ref{fig:noise}.

 

\begin{figure}[ht]
  \begin{center}
 \includegraphics[scale=0.3]{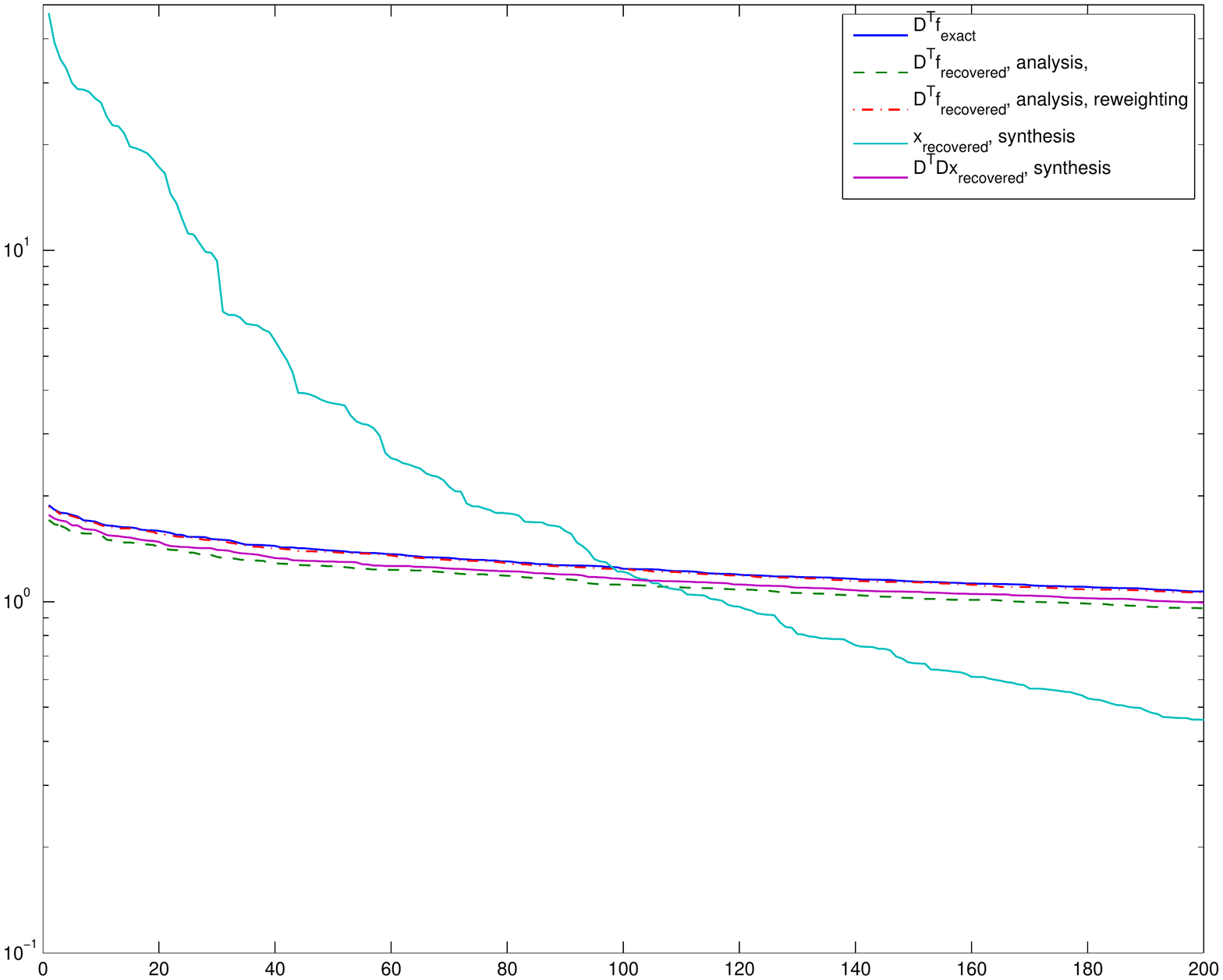}
\end{center}
\caption{The largest $200$ coefficients of the coefficient vector
  $D^*f$ (blue), $D^*\hat{f}$ from $\ell_1$-analysis (dashed green),
  $D^*\hat{f}$ from $\ell_1$-analysis with $3$ reweighting iterations
  (dashed red), $\hat x$ from $\ell_1$-synthesis (cyan), and
  $D^*\hat{f}$ from $\ell_1$-synthesis (magenta).}\label{fig:SynAn}
\end{figure}

An alternative to $\ell_1$-analysis is $\ell_1$-synthesis, which we
discuss in Section~\ref{sec:alts}; $\ell_1$-synthesis minimizes in the
coefficient domain, so its solution is a vector $\hat{x}$, and we set
$\hat{f} = D\hat{x}$.  Our next simulation confirms that although we
cannot recover the coefficient vector $x$, we can still recover the
signal of interest.  Figure~\ref{fig:SynAn} shows the largest $200$
coefficients of the coefficient vector $x$, and those of $D^*f$ as
well as $D^*\hat{f}$ for both $\ell_1$-analysis and
$\ell_1$-synthesis.  The plot also shows that the recovery of
$\ell_1$-analysis with reweighting outperforms both standard
$\ell_1$-analysis and $\ell_1$-synthesis.

\begin{figure}[ht]
\begin{center}
  \includegraphics[scale=0.4]{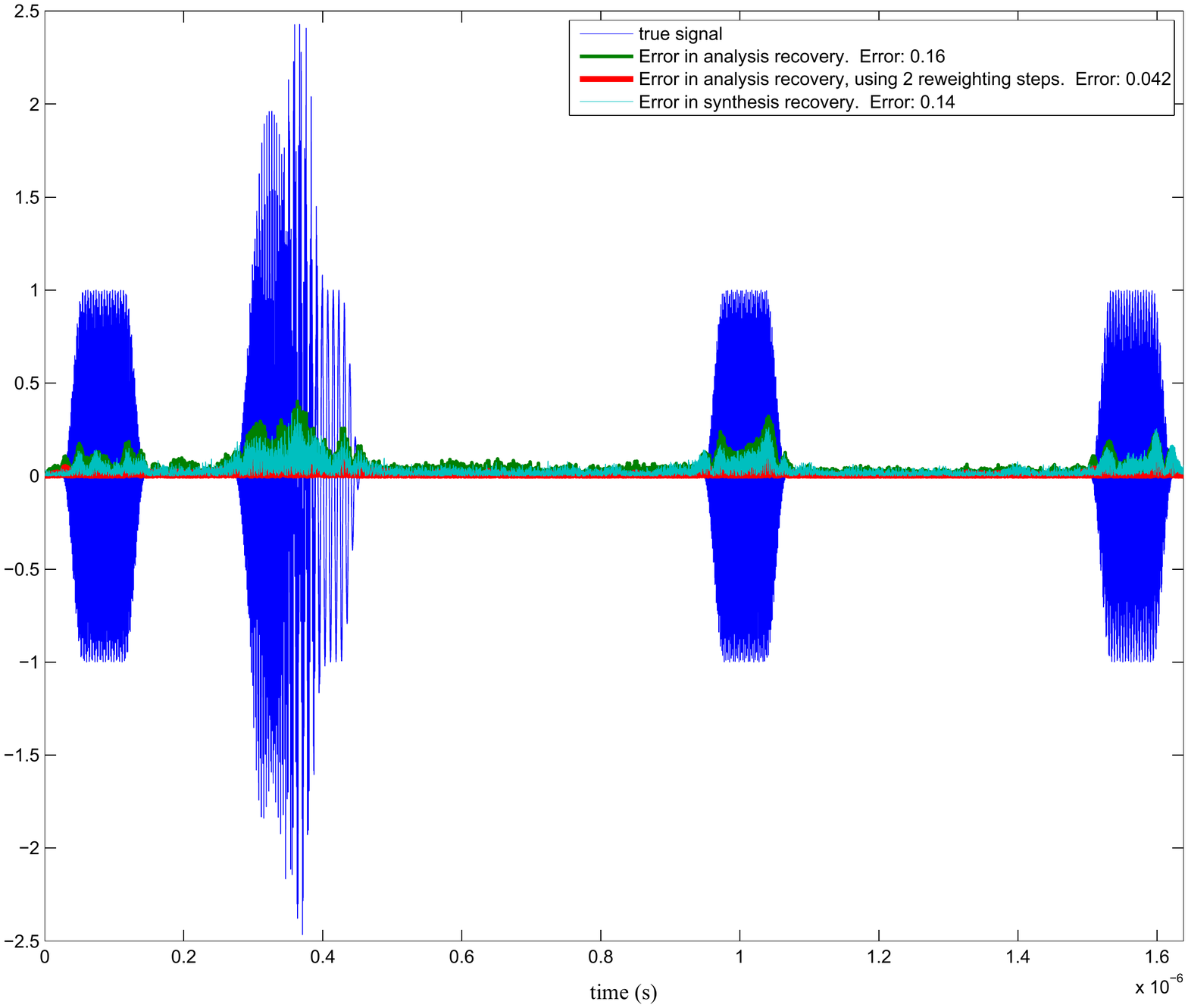}
\end{center}
\caption{Recovery (without noise) of a compressible signal in the time
  domain.  Blue denotes the actual signal, while green, red, and cyan
  denote the recovery error from $\ell_1$-analysis, reweighted
  $\ell_1$-analysis ($2$ iterations), and $\ell_1$-synthesis,
  respectively. The legend shows the relative error $\|\hat f -
  f\|_2/\|f\|_2$ of the three methods.}\label{fig:comp}
\end{figure}

Our final simulation compares recovery error on a compressible signal
(in the time domain) for the $\ell_1$-analysis, reweighted
$\ell_1$-analysis, and $\ell_1$-synthesis methods.  We see in
Figure~\ref{fig:comp} that the $\ell_1$-analysis and
$\ell_1$-synthesis methods both provide very good results, and that
reweighted $\ell_1$-analysis provides even better recovery error.

\section{Discussion}\label{sec:disc}

Theorem~\ref{thm:main} shows that $\ell_1$-analysis is accurate when
the coefficients of $D^*f$ are sparse or decay rapidly.  As discussed
above, this occurs in many important applications.  However, if it is
not the case, then the theorem does not guarantee good recovery.  As
previously mentioned, this may occur when the dictionary $D$ is a
concatenation of two (even orthonormal) bases.  For example, a signal
$f$ may be decomposed as $f = f_1 + f_2$ where $f_1$ is sparse in the
basis $D_1$ and $f_2$ is sparse in a different basis, $D_2$.  One can
consider the case where these bases are the coordinate and Fourier
bases, or the curvelet and wavelet bases, for example.  In these
cases, $D^*f$ is likely to decay slowly since the component that is
sparse in one basis is not at all sparse in the
other~\cite{DK10:Microlocal}.  This suggests that $\ell_1$-analysis
may then not be the right algorithm for reconstruction in such
situations.

\subsection{Alternatives}\label{sec:alts}
Even though $\ell_1$-analysis may not work well in this type of setup,
one should still be able to take advantage of the sparsity in the
problem.  We therefore suggest a modification of $\ell_1$-analysis
which we call \textit{Split-analysis}.  As the name suggests, this
problem splits up the signal into the components we expect to be
sparse:
\begin{equation*}
  (\hat{f}_1, \hat{f}_2) = \argmin_{\tilde{f}_1, \tilde{f}_2} \|D_1^*\tilde{f}_1\|_1 + \|D_2^*\tilde{f}_2\|_1 \quad\text{subject to}\quad \|A(\tilde{f}_1 + \tilde{f}_2) - y\|_2 \leq \eps.
\end{equation*}
The reconstructed signal would then be $\hat{f} = \hat{f}_1 +
\hat{f}_2$.  Some applications of this problem in the area of image restoration have been
studied in~\cite{COS09:Split}.  Since this is an analagous problem to $\ell_1$-analysis,
one would hope to have a result for Split-analysis similar to
Theorem~\ref{thm:main}.

An alternative way to exploit the sparsity in $f = f_1 + f_2$ is to
observe that there may still exist a (nearly) sparse expansion $f = Dx
= D_1 x_1 + D_2 x_2$.  Thus one may ask that if the coefficient vector
$x$ is assumed sparse, why not just minimize in this domain?  This
reasoning leads to an additional approach, called $\ell_1$-Synthesis
or Basis Pursuit (see also the discussion in~\cite{EMR07:Analysis}):
\begin{equation}\tag{$\ell_1$-synthesis}
  \hat{x} = \argmin_{\tilde{x}} \|\tilde x\|_1 \quad\text{subject to}\quad \|AD\tilde x - y \|_2 \leq \eps.
\end{equation} 
The reconstructed signal is then $\hat{f} = D\hat{x}$.  Empirical studies also show that $\ell_1$-synthesis often provides good recovery, however, it is fundamentally distinct from $\ell_1$-analysis.  The geometry of the two problems is analyzed in~\cite{EMR07:Analysis}, and there it is shown that because these geometrical structures exhibit substantially different properties, there is a large gap between the two formulations.  This theoretical gap is also demonstrated by numerical simulations in~\cite{EMR07:Analysis}, which show that the two methods perform very differently on large families of signals.

\subsection{Fast Transforms}\label{sec:Fast}
For practical reasons, it is clearly advantageous to be able to use
measurement matrices $A$ which allow for easy storage and fast
multiplication.  The partial DFT for example, exploits the Fast
Fourier Transform (FFT) which allows the sampling matrix to be
applied to a $n$-dimensional vector in $\bigO(n \log n)$ time, and
requires only $\bigO(m\log n)$ storage.  Since the partial
DFT has been proven to satisfy the RIP~\cite{CT04:Near-Optimal} (see also~\cite{RV08:sparse}), it is a
fast measurement matrix that can be used in many standard compressed
sensing techniques.  

One of course hopes that fast measurement
matrices can be used in the case of redundant and coherent
dictionaries as well.  As mentioned, the result of Krahmer and Ward implies
that any matrix which satisfies the RIP will satisfy the D-RIP when multiplied
by a random sign matrix~\cite{KW10:NewAnd}.  Therefore, the $m\times n$ subsampled Fourier matrix with
$m = \bigO(s\log^4n)$ along
with the sign matrix will satisfy the D-RIP and provides the fast multiply.

A result at the origin of this notion was proved by Ailon and Liberty, also
after our initial submission of this paper~\cite{AL10:AnAlmost}.  Recall that in light of~\eqref{eq:cond}, we desire
a fast transform that satisfies the Johnson-Lindenstrauss lemma in the 
following sense.  For a set $Q$ of $N$ vectors in $n$-dimensional space,
we would like a fast transform $A$ that maps this set into a space of dimension $\bigO(\log N)$ 
(possibly with other factors logarithmic in $n$) such that
\begin{equation}\label{JLwant}
(1-\delta)\|v\|_2^2 \leq \|Av\|_2^2 \leq (1+\delta)\|v\|_2^2 \quad\text{ for all $v\in Q$.}
\end{equation}
Note that the dimension $\bigO(\log N)$ will of course also depend on the constant $\delta$.  Due to standard covering 
arguments (see e.g.~\cite{RSV08:Redundant,BDDW07:Johnson-Lindenstrauss}), this would yield 
a $m \times n$ fast transform with optimal
number of measurements, $m = \bigO(s\log n)$, obeying the D-RIP.   

 Ailon and Liberty show that the subsampled Fourier matrix multiplied by a random sign matrix does exactly this~\cite{AL10:AnAlmost}.
Thus in other words,
for a fixed $m$, this $m \times n$ construction satisfies the D-RIP up to sparsity level $s = \bigO(m/\log^4 n)$.
The cost of a matrix--vector
multiply is of course dominated by that of the FFT, $\bigO(n\log n)$. 
Its storage requirements are also $\bigO(m\log n)$. 
Their results can also be generalized to other transforms with the same type of fast multiply.  

These results yield a transform with a fast multiply which satisfies the D-RIP.  The number of measurements
and the multiply and storage costs of the matrix are of the same magnitude as those that satisfy the RIP.
The D-RIP is, therefore, satisfied by matrices with the same benefits as those in standard
compressed sensing.  This shows that compressed sensing with redundant and coherent
dictionaries is viable with completely the same advantages as in the standard setting.


\subsection*{Acknowledgements}

This work is partially supported by the ONR grants N00014-10-1-0599
and N00014-08-1-0749, the Waterman Award from NSF, and the NSF DMS
EMSW21-VIGRE grant. EJC would like to thank Stephen Becker for
valuable help with the simulations.

\small 
\bibliography{thebib}

\end{document}